\documentclass[3p]{elsarticle}

\usepackage{amsmath} 
\usepackage{amssymb}
\usepackage{amsthm}
\usepackage{amsfonts}
\usepackage{bm}
\graphicspath{ {./figs/} }

\newproof{pf}{Proof}

\numberwithin{equation}{section}
\journal{arXiv}

\begin{document}

\begin{frontmatter}

\title{Nonlinear approximation of functions based on non-negative least squares solver}

\author[nsi,uni]{Petr N. Vabishchevich\corref{cor}}
\ead{vabishchevich@gmail.com}
\ead[url]{https://sites.google.com/view/vabishchevich/}

\address[nsi]{Nuclear Safety Institute, Russian Academy of Sciences,
              52, B. Tulskaya, 115191 Moscow, Russia}
\address[uni]{North-Caucasus Center for Mathematical Research, North-Caucasus Federal University, 
                   1, Pushkin Street, 355017 Stavropol, Russia}

\cortext[cor]{Corresponding author}

\begin{abstract}
In computational practice, most attention is paid to rational approximations of functions and approximations by the sum of exponents.
We consider a wide enough class of nonlinear approximations characterized by a set of two required parameters. The approximating function is linear in the first parameter; these parameters are assumed to be positive.
The individual terms of the approximating function represent a fixed function that depends nonlinearly on the second parameter.
A numerical approximation minimizes the residual functional by approximating function values at individual points.
The second parameter's value is set on a more extensive set of points of the interval of permissible values.
The proposed approach's key feature consists in determining the first parameter on each separate iteration of the classical non-negative least squares method.
The computational algorithm is used to rational approximate the function $x^{-\alpha}, \ 0 < \alpha < 1, \ x \geq 1$.
The second example concerns the approximation of the stretching exponential function $\exp(- x^{\alpha} ), \ \ \quad 0 < \alpha < 1$ at $ x \geq 0$ by the sum of exponents. 
\end{abstract}



\end{frontmatter}

\section{Introduction}

Problems of nonlinear approximation of functions are in high demand in computational practice for many applied problems. Rational approximations are widely used in various variants \cite{braess1986nonlinear}.
Remez algorithm \cite{iske2018approximation} and recently developed algorithms for rational approximation \cite{nakatsukasa2018aaa,hofreither2021algorithm} are used to find parameters of the approximating function.
Much attention is paid to approximation by the sum of exponents in function approximation.
In particular, achievements in this area are reflected by the works of \cite{holmstrom2002review,hristov2022prony}.
In many cases, the approximation of the function $f(x)$ is
\[
f(x) \approx \sum_{i=1}^{m} u_i \varphi(x, v_i)
\]
with the known functional dependence $\varphi(x, v)$.
Some restrictions may be imposed on the approximation parameters to be sought.
A typical example with non-negative coefficients $u_i, \ i = 1,2, \dots, m$ is interesting for many applications.

The theory \cite{mhaskar2000fundamentals} and computational practice \cite{meinardus2012approximation,trefethen2019approximation} of function approximation are well developed in the linear case.
In this the function $f(x)$ is approximated by a given set of trial functions $\varphi_i(x), \ i = 1,2, \dots, m$:
\[
f(x) \approx \sum_{i=1}^{m} u_i \varphi_i(x) .
\]
Optimal approximations are constructed in Hilbert spaces using the least squares method \cite{lawson1995solving,bjorck1996numerical}.
Special note that computational algorithms have long been well developed to consider the constraints $u_i > 0, \ i = 1,2, \dots, m$.
We propose a heuristic algorithm for solving the nonlinear function approximation problem.
It is based on expanding the set of test functions and their subsequent selection during the iterative solution of the nonlinear least squares problem.

The paper is organized as follows.
In Section 2, the problem of nonlinear approximation of functions is posed.
The proposed computational algorithm is described in Section 3.
Section 4 presents results on rational approximation and approximation by the sum of exponents of functions $x^{-\alpha}, \ x \geq 1$ and $\exp(- x^{\alpha} ), \ x \geq 0$ at $0 < \alpha < 1$.
The work results are summarized in Section 5.

\section{Problem formulation}

We consider the problem of nonlinear approximation of one-dimensional function $f(x)$ on the interval $[a, b]$.
In the Hilbert space $L_2([a,b], \varrho (x))$ with weight $\varrho (x) > 0$ the scalar product and norm are defined as follows
\[
 (g, q) = \int_{a}^{b} \varrho (x) \big (g(x) - q(x) \big )^2 d x ,
 \quad \| g \| = (g , g )^{1/2} .
\] 
The function $f(x)$ is approximated by the function $r(x, \bm u, \bm v)$ with two numerical parameter sets $\bm u = \{u_1, u_2, \ldots, u_m\}$, $\bm v = \{v_1, v_2, \ldots, v_m\}$.
Let's assume that the approximating function has the form
\begin{equation}\label{1}
 r(x, \bm u, \bm v) = \sum_{i=1}^{m} u_i \varphi(x, v_i) .
\end{equation} 
In the representation (\ref{1}) we have isolated the linear coefficients $u_i, \ i =1,2, \ldots, m$, and the dependence on $v_i, \ i =1,2, \ldots, m,$ is determined by the given function $\varphi(x,v)$.

Among the methods of nonlinear approximation of functions, rational approximation and approximation by the sum of exponents are the most widespread.
In particular, in the case of rational approximation at $a \geq 0$ we use the parametric function 
\begin{equation}\label{2}
 \varphi(x, v) = \frac{1}{1 + v x} .
\end{equation} 
When approximated by the sum of the exponents, we have
\begin{equation}\label{3}
 \varphi(x, v) = \exp(- v x) .
\end{equation} 

In the nonlinear approximation problem we consider, $\bm u, \bm v$ are subject to the following restrictions: the parameters $u_i, \ i =1,2, \ldots, m,$ are non-negative, and $v_i, \ i =1,2, \ldots, m,$ are chosen from the interval $[c, d]$. We come to the problem
\begin{equation}\label{4}
 J(\bm u, \bm v) \rightarrow \min, 
 \quad   (\bm u, \bm v) \in K,
\end{equation} 
where 
\begin{equation}\label{5}
 J(\bm u, \bm v) = \|f(x) - r(x, \bm u, \bm v)\|^2,
 \quad K = \{ (\bm u, \bm v) \ | \ u_i > 0, \ v_i \in [c, d], \ i =1,2, \ldots, m \} .
\end{equation} 
An approximate solution to the functional minimization problem (\ref{1}), (\ref{4}), (\ref{5}) is constructed by setting the functions $f(x), \ r(x, \bm u, \bm v)$ on a sufficiently detailed set of points on the interval $[a,b]$. 

\section{Nonlinear approximation algorithm}

We begin by dividing the interval $[a,b]$ into $n$ partial intervals of length $h_j, \ j = 1,2, \ldots, n$, so that
\[
 b - a = \sum_{j=1}^{n} h_j.
\] 
We denote the centers of the intervals $h_j$ by $x_j, \ j = 1,2, \ldots, n$.
Using the quadrature formula of rectangles, we compare the original problem (\ref{4}), (\ref{5}) with the minimization problem
\begin{equation}\label{6}
 J^h(\bm u, \bm v) \rightarrow \min, 
 \quad   (\bm u, \bm v) \in K,
\end{equation} 
\begin{equation}\label{7}
 J^h(\bm u, \bm v) = \sum_{j=1}^{n} \varrho(x_j) \big (f(x_j) - r(x_j, \bm u, \bm v) \big )^2 h_j .
\end{equation} 

To simplify the approximation problem, the coefficients $v_i, \ i =1,2, \ldots, m,$ will not be evaluated over the whole interval $[c, d]$, but only at given points
\[
 \widetilde{v}_k \in V_l,
 \quad V_l = \{c \leq \widetilde{v_1} \leq \widetilde{v_2} \leq \ldots \leq \widetilde{v_l} \leq d \}
\]
for a sufficiently fine partition ($l \gg m$). 
Thus we take instead of $K$ the set of constraints in the form
\[
 \widetilde{K} = \{ (\bm u, \bm v) \ | \ u_i > 0, \ v_i \in V_l, \ i =1,2, \ldots, m \} .
\] 
Given (\ref{1}), we will come from the problem (\ref{6}), (\ref{7}) to the problem 
\begin{equation}\label{8}
 J^h(\bm u, \bm v) \rightarrow \min, 
 \quad   (\bm u, \bm v) \in \widetilde{K},
\end{equation} 
\begin{equation}\label{9}
J^h(\bm u, \bm v) =  \sum_{j=1}^{n} \varrho(x_j)  h_j \Big ( \sum_{i=1}^{m} u_i \varphi(x_j, v_i) -  f(x_j) \Big )^2 .
\end{equation} 

Instead of the vector $\bm u$ with components $u_i, \ i = 1,2, \ldots, m,$ we introduce a vector $\widetilde{\bm u}$ of larger dimension with components $\widetilde{u}_k, \ k = 1,2, \ldots, l$. The components of the vector $\widetilde{\bm u}$ we define from the approximation condition 
\begin{equation}\label{10}
  \sum_{i=1}^{m} u_i \varphi(x_j, v_i) = \sum_{k=1}^{l} \widetilde{u}_k \varphi(x_j, \widetilde{v}_k) . 
\end{equation} 
By doing so, we put $\widetilde{u}_k = u_i$ if $\widetilde{v}_k = v_i$ and $\widetilde{u}_k = 0$ if $\widetilde{v}_k \neq v_i$.
Considering (\ref{10}), we proceed from the minimization problem (\ref{8}), (\ref{9}) to the problem
\begin{equation}\label{11}
 \widetilde{J} (\widetilde{\bm u}) \rightarrow \min, 
 \quad  \widetilde{u}_k \geq 0, 
 \quad k = 1,2, \ldots, l ,
\end{equation} 
\begin{equation}\label{12}
\widetilde{J} (\widetilde{\bm u})  =  \sum_{j=1}^{n} \varrho(x_j)  h_j \Big ( \sum_{k=1}^{l} \widetilde{u}_k \varphi(x_j, \widetilde{v}_k) -  f(x_j) \Big )^2 .
\end{equation} 
As a result, we obtained a constrained (non-negative) least squares problem of a larger dimension to determine the parameters $\widetilde{u}_k, \ k = 1,2, \ldots, l,$ in the linear representation of the approximating function.

Computational algorithms for the minimization problem (\ref{11}), (\ref{12}) are well studied \cite{bjorck1996numerical,lawson1995solving}.
In computational practice, the most widely used is the NNLS (Non-Negative Least Squares) algorithm described in detail in \cite{lawson1995solving}.
Given the specifics of the problem (\ref{11}), (\ref{12}), we will separately specify variants of NNLS algorithms for large-scale problems (see, for example, \cite{kim2013non,myre2017tnt}).

The standard NNLS algorithm is a two-step iterative method with main and inner loop iterations.
The number of positive coefficients is initially set to zero.
As the number of iterations increases, the residual decreases, though not monotonically.
First, decreasing the residual is provided by increasing the number of positive coefficients.

Taking into account the above features of the iterative process of NNLS algorithm, we can propose the following strategy of coefficient selection in a nonlinear approximation of functions based on the minimization problem (\ref{11}), (\ref{12}).
We perform a sufficiently large number of iterations of the NNLS algorithm.
At each iteration, we control the residual and the number $\widetilde{m}$ of positive coefficients $\widetilde{u}_k, \ k = 1,2, \ldots, l$.
The number of iterations of the NNLS algorithm is chosen such that $\widetilde{m} = m$ and the residual is minimal.
No additional modifications of the standard algorithm are performed.
The computational implementation is based on the non-negative least squares solver from the SciPy library \cite{2020SciPy-NMeth} (module \textsf{optimize}, function \textsf{nnls}).

\section{Numerical experiments} 

Let us illustrate the possibility of constructing nonlinear approximations of functions with two examples.
First, we construct rational approximations of the function $x^{-\alpha}, 0 < \alpha < 1$ at $x \geq 1$.
In the last decade, such a problem has been actively discussed in the literature in connection with solving boundary value problems with fractional power elliptic operators (see, for example, \cite{bonito2018numerical,harizanov2020rev}), and also when considering more general problems with \cite{vabishchevich2022some} operator functions.
The second example concerns the approximation of the function $\exp(-x^{-\alpha}), 0 < \alpha < 1$ at $x \geq 0$.
In this case, we use approximations by the sum of exponents.
Such problems are typical in approximate solutions of nonstationary problems with memory when the difference kernel of the integral term \cite{vabMemory} is approximated.

\subsection{Approximation of $x^{-\alpha}$} 

We will approximate the function $f(x) = x^{-\alpha}$ when $a = 1$ and $b$ is large enough.
We used $b = 10^{15}$ in the following calculations.
In applied problems, it is often important to impose an additional restriction on the approximating function
\begin{equation}\label{13}
r(a, \bm u, \bm v) = f(a) .
\end{equation}
In the class of rational approximations of the type (\ref{2}), given (\ref{13}), for the approximation $x^{-\alpha}$, we obtain the representation
\[
x^{-\alpha} \approx 1 + \sum_{i=1}^{m} u_i \varphi(x, v_i) ,
\]
when
\[
 \varphi(x, v_i) = \frac{1}{1 + v_i x} - \frac{1}{1 + v_i},
 \quad i = 1,2, \ldots, m . 
\] 

The choice of the weight function $\varrho (x)$, which allows controlling the approximation accuracy in some parts of the interval $[a,b]$, requires special attention.
In our case, the approximated function decreases to zero as $x$ increases.
A more significant influence of points at small $x$ is provided by setting a decreasing function $\varrho (x)$.
To partition the interval $[a,b]$ we use partial intervals of increasing length ($h_{i+1} > h_i, \ i = 1,2, \ldots, m-1$).
This approach is formalized by introducing a new variable $\theta$ instead of $x$.
Put $x = \exp(\theta)$, so that $\theta \in [0, \beta]$ ($\exp(\beta) = b$).
For the residual functional, we get
\[
J(\bm u, \bm v) = \int_{0}^{\beta }\varrho \big(\exp(\theta) \big) \exp(\theta) \Big (\exp(-\alpha\theta) - 1
- \sum_{i=1}^{m} u_i \varphi \big(\exp(\theta), v_i \big) \Big )^2 d \theta .
\]
The following is the computational data when setting $\varrho \big(\exp(\theta) \big) = \exp(-\theta)$ ($\varrho(x) = x^{-1}$).

\begin{figure}[htbp]
\centering
\includegraphics[width=0.49\linewidth]{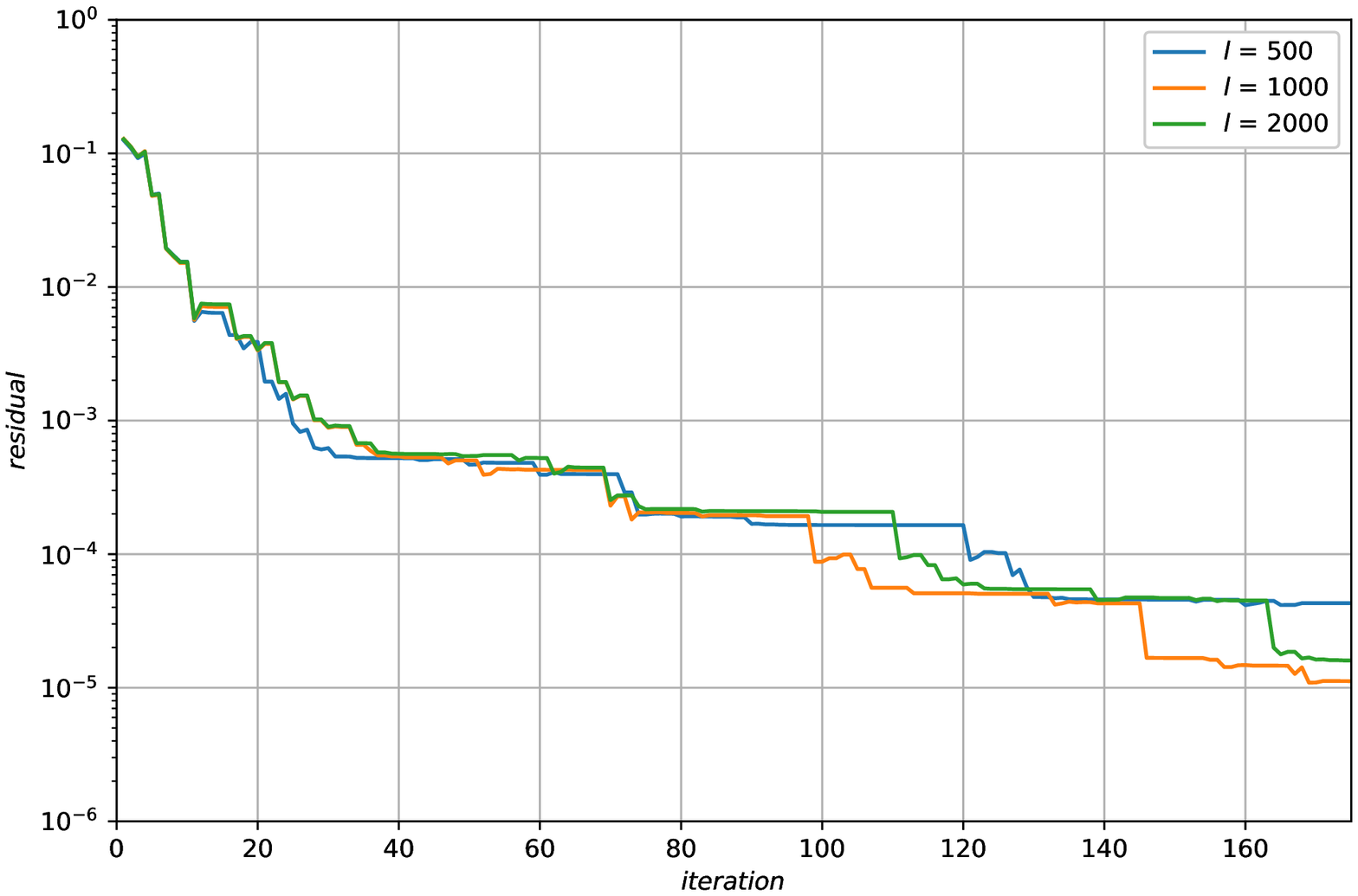} 
\includegraphics[width=0.49\linewidth]{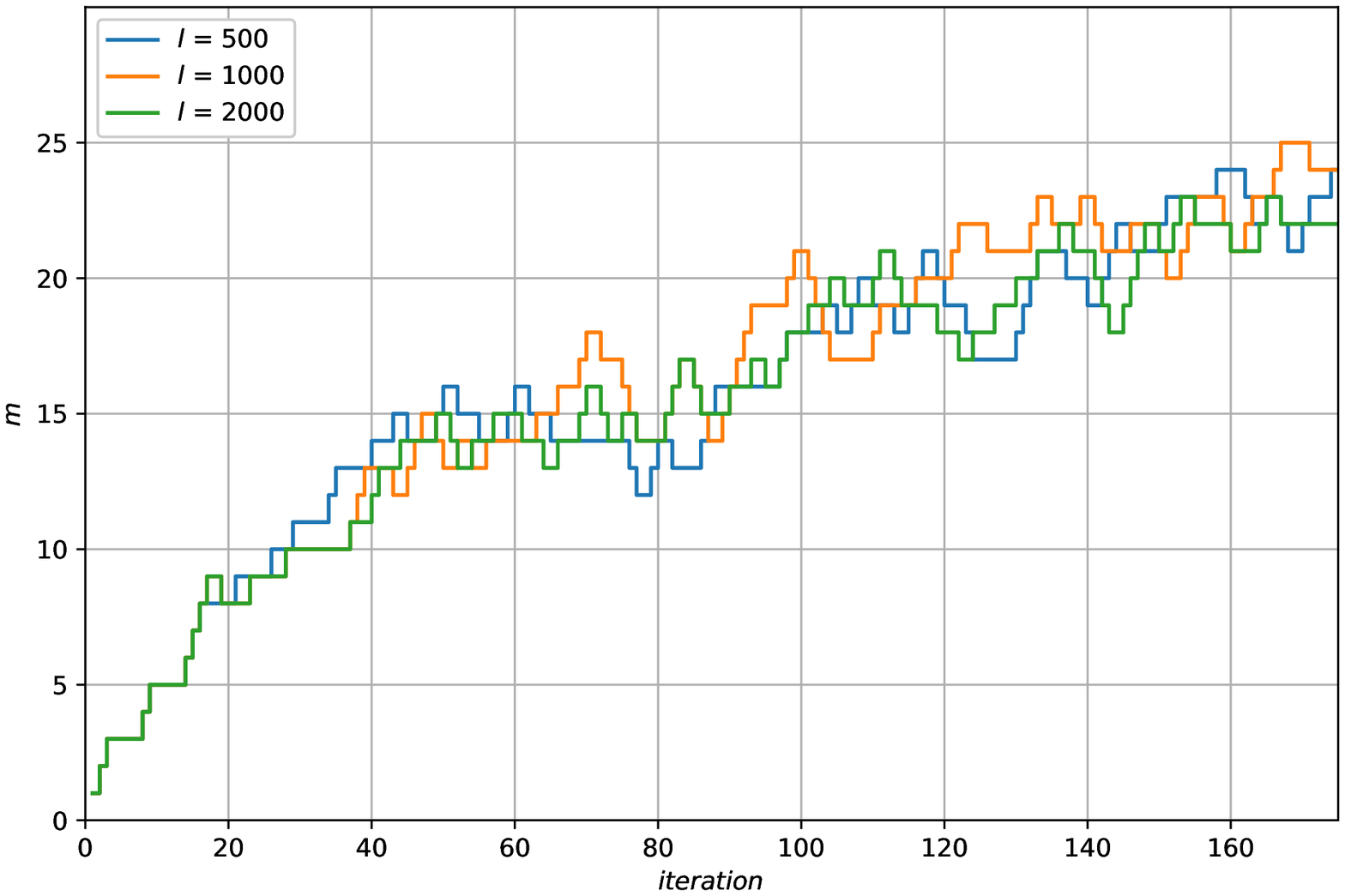} 
\caption{Residual (left) and the number of non-zero elements $m$ (right) for $\widetilde{u}_k, \ k = 1,2, \ldots, l,$ in individual iterations of the non-negative least squares method.}
\label{fig-1}
\end{figure}

At sufficiently large values of the number of points $n$ on the interval $[a,b]$, the approximation accuracy changes insignificantly.
In the presented results of calculations, we limited ourselves to the case $n=5000$.
Of greater interest are the calculated data at different partitioning of the interval $[c,d]$ of permissible values of the parameter $v$.
Figure~\ref{fig-1} illustrates the effect of $l$ in (\ref{12}) when we approximate the function $x^{-\alpha}$ with $\alpha = 0.5$.
The residual equal to $\widetilde{J}^{1/2} (\widetilde{\bm u})$ using different numbers of iterations of the non-negative least squares method is shown in the left-hand side of the figure. We observe a reasonably fast, generally speaking, non-monotonic decrease in the residual with the increasing number of iterations.
This residual is achieved with varying numbers of non-zero coefficients $m$ with a general tendency of $m$ increasing as the number of iterations of the non-negative least squares method increases. From these data, we estimate the number of iterations to achieve the minimum norm of residual for a given number $m$ of terms of the rational approximation (\ref{10}).

\begin{figure}[htbp]
\centering
\includegraphics[width=0.49\linewidth]{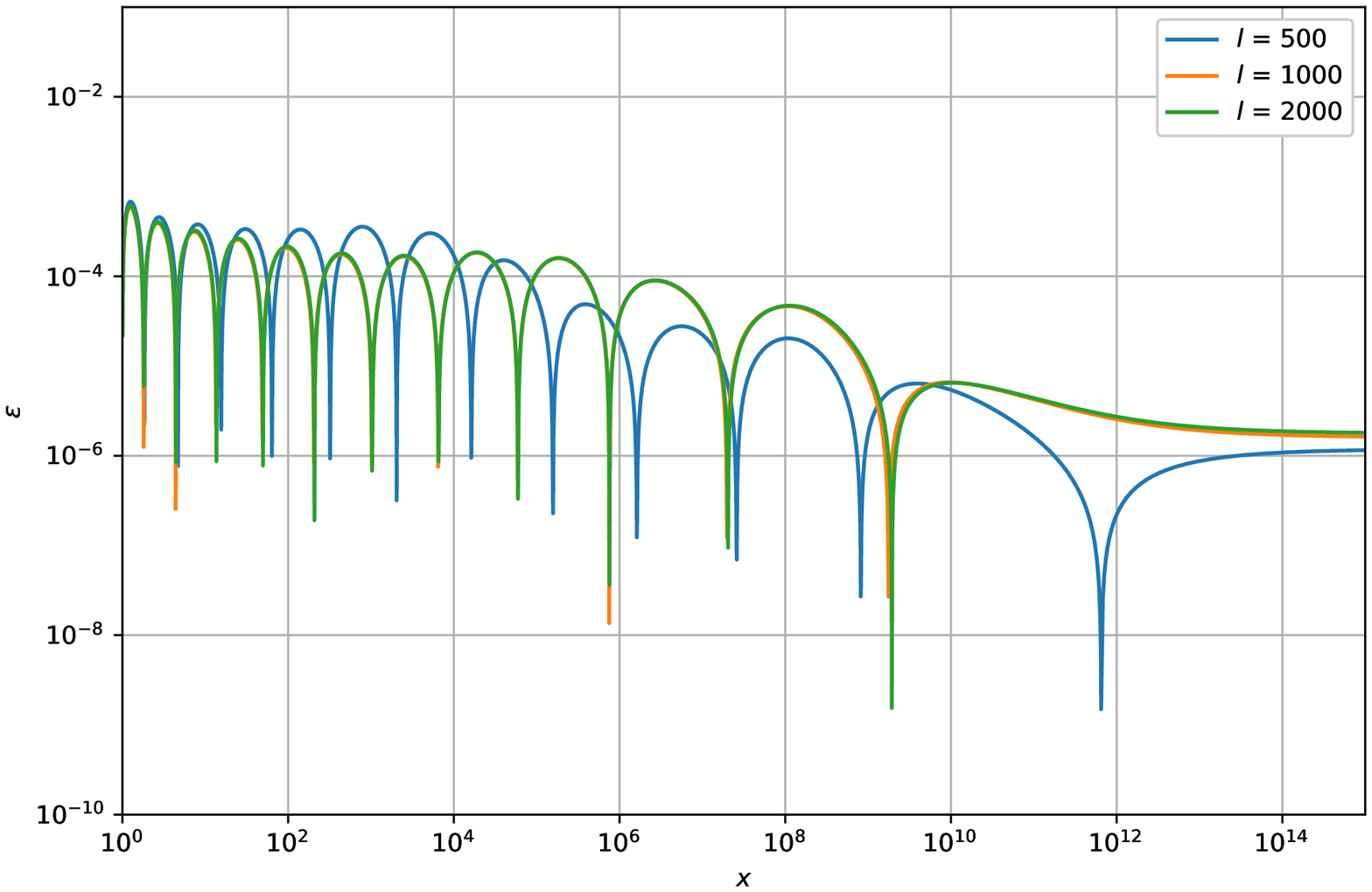} 
\includegraphics[width=0.49\linewidth]{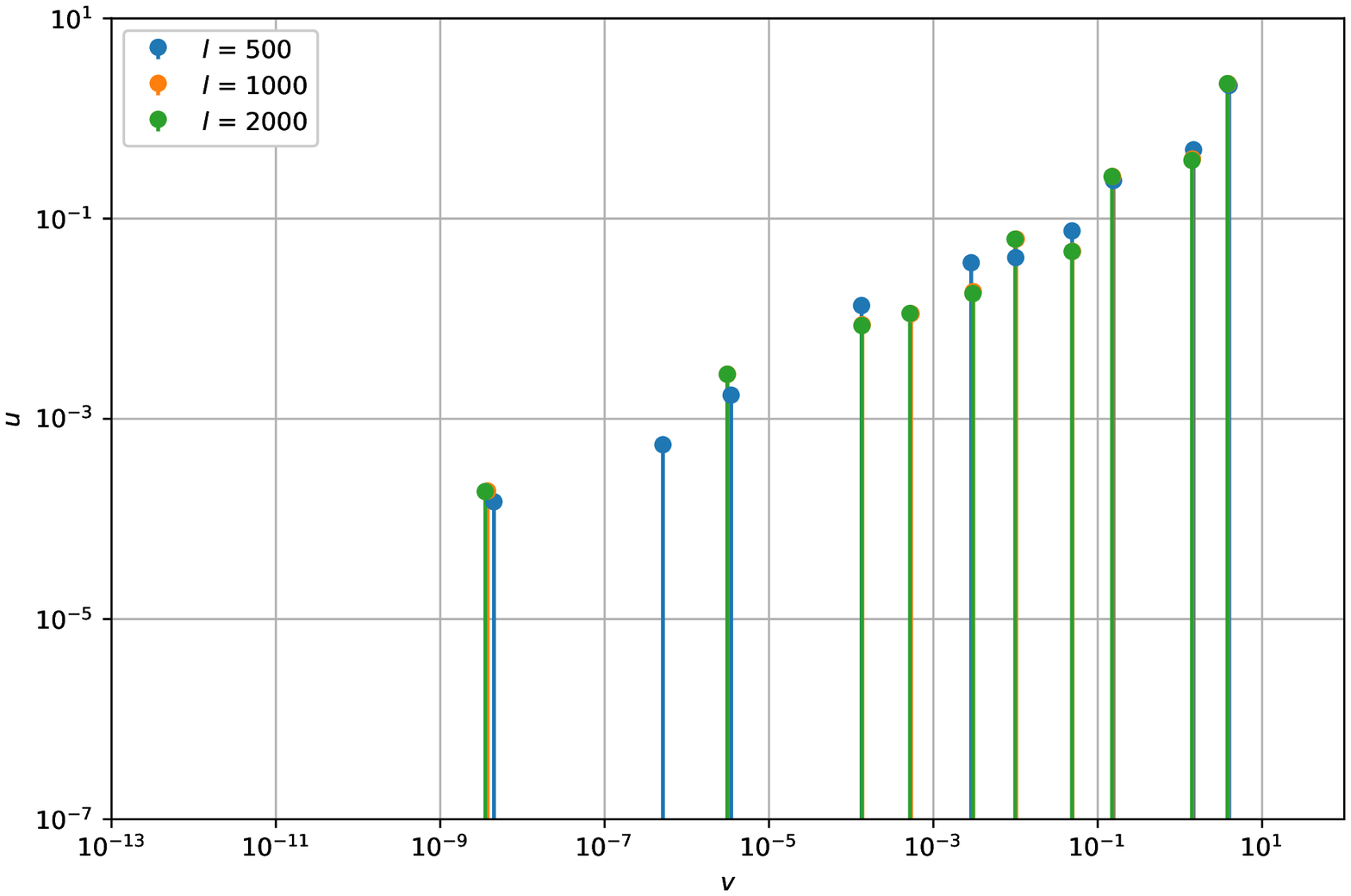} 
\caption{Approximation accuracy $\varepsilon$ (left) and approximation parameters $u_i, v_i, \ i = 1,2, \ldots, m,$ (right) at different partitioning of the interval $[c,d]$.}
\label{fig-2}
\end{figure}

The effect of interval partitioning $[c,d]$ on the approximation accuracy is shown in Fig.\ref{fig-2}.
Here is the data for $m=10$ at $l = 500, 1000, 2000$.
The approximation accuracy is estimated by the value
\[
\varepsilon (x_j) = \Big | \sum_{i=1}^{m} u_i \varphi(x_j, v_i) - f(x_j) \Big |,
\quad j = 1,2, \ldots, n .
\]
We observe a remarkable similarity in determining the approximation parameters $u_i, v_i, \ i = 1,2, \ldots, m,$ (right-hand side of the figure) and in the accuracy of the approximation function $x^{-\alpha}$ with $\alpha=0.5$ at $x \in [1, 10^{15}]$ (left-hand side of the figure).
With this in mind, we fixed the number of partitions $l=1000$ in our calculations.

\begin{figure}[htbp]
\centering
\includegraphics[width=0.49\linewidth]{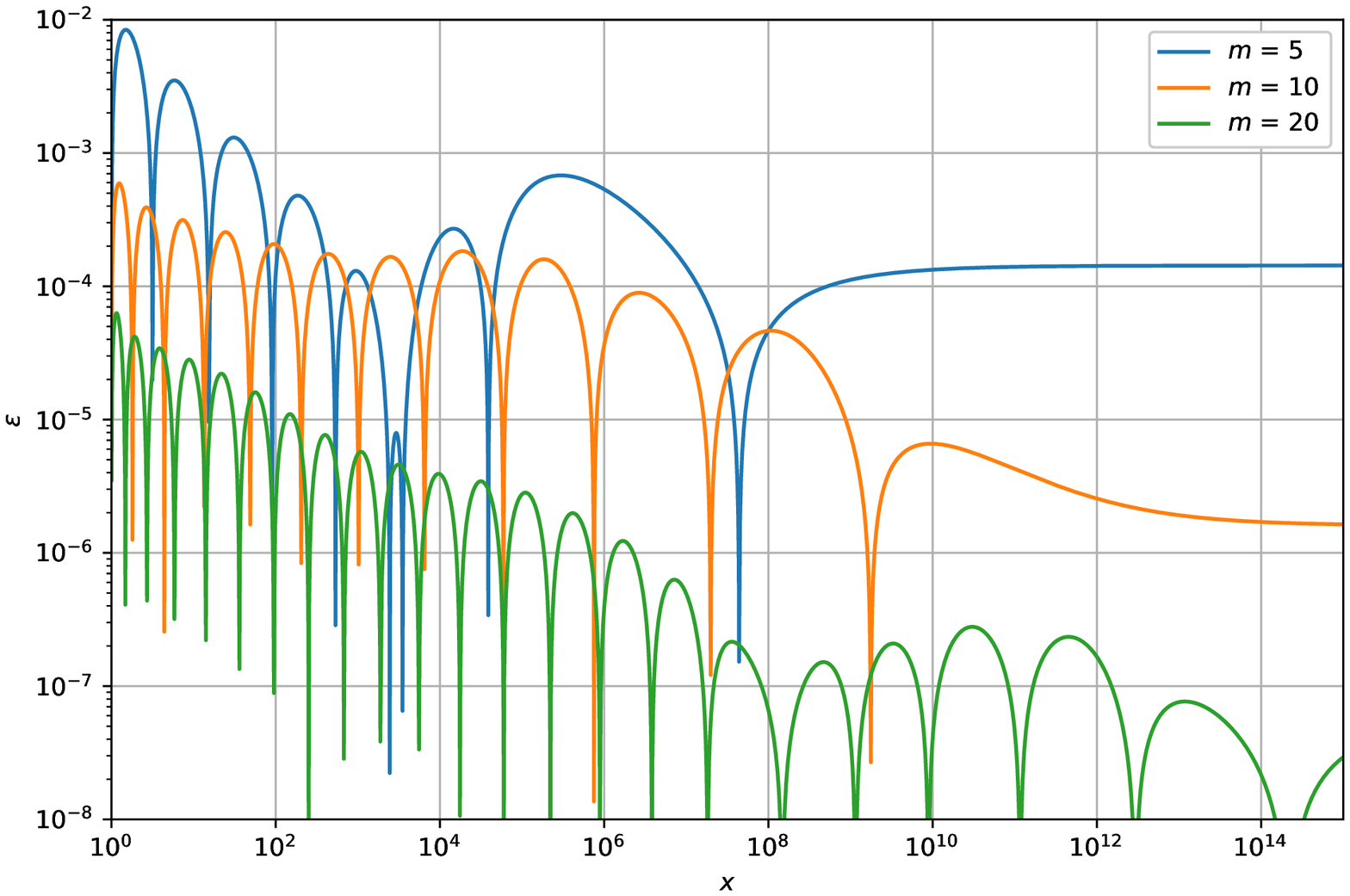} 
\includegraphics[width=0.49\linewidth]{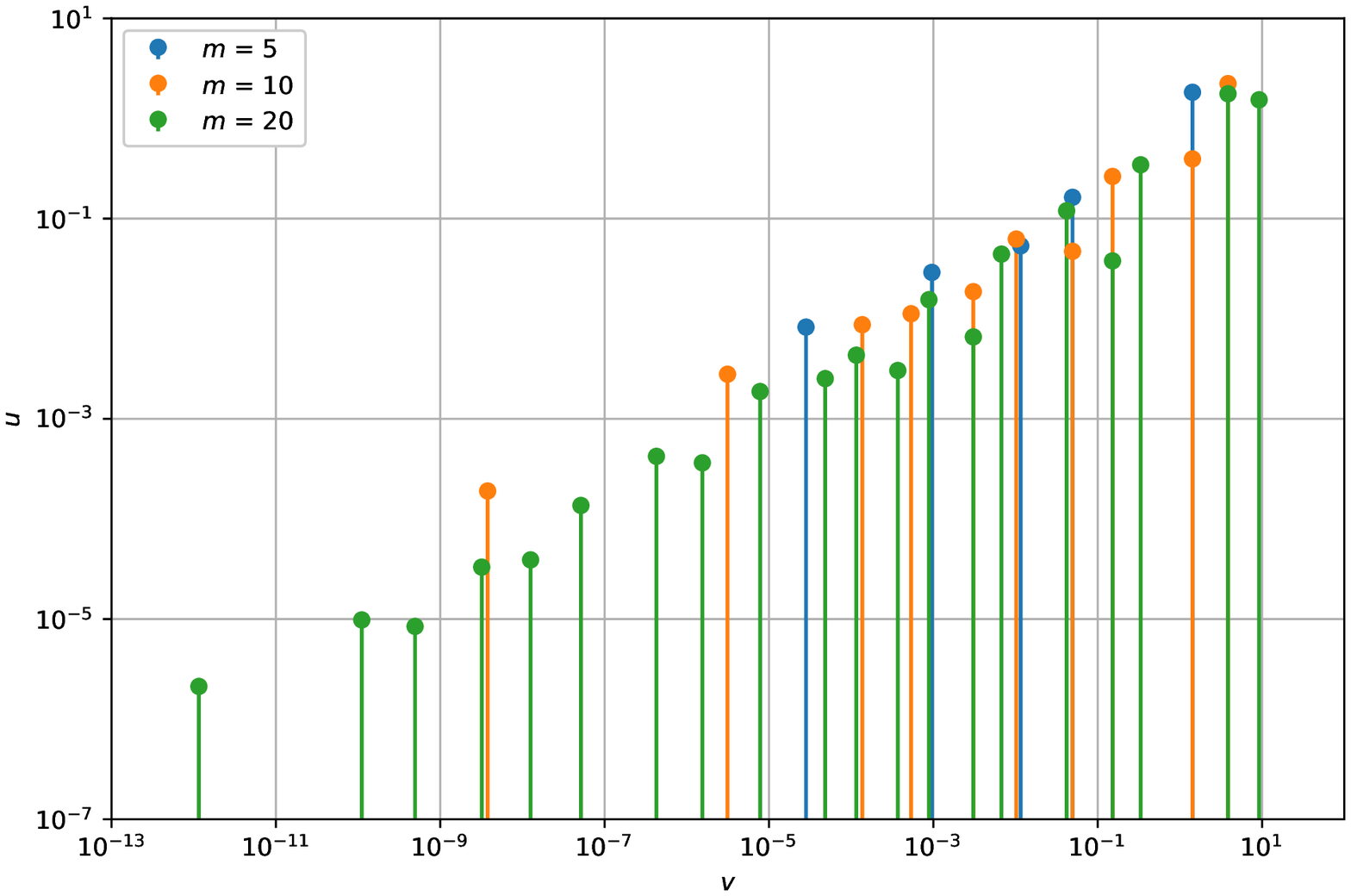} 
\caption{Approximation accuracy $\varepsilon$ (left) and approximation parameters $u_i, v_i, \ i = 1,2, \ldots, m,$ (right) for the function $x^{-\alpha}$ with $\alpha = 0.5$ at various $m$.} 
\label{fig-3}
\end{figure}

Figure~\ref{fig-3} shows the approximation parameters $u_i, v_i, \ i = 1,2, \ldots, m,$ when given $m = 5, 10, 20$.
There is a significant increase in the approximation accuracy with increasing $m$.
For this variant, the number of iterations is $13, 38$, and $110$, respectively.

\begin{figure}[htbp]
\centering
\includegraphics[width=0.49\linewidth]{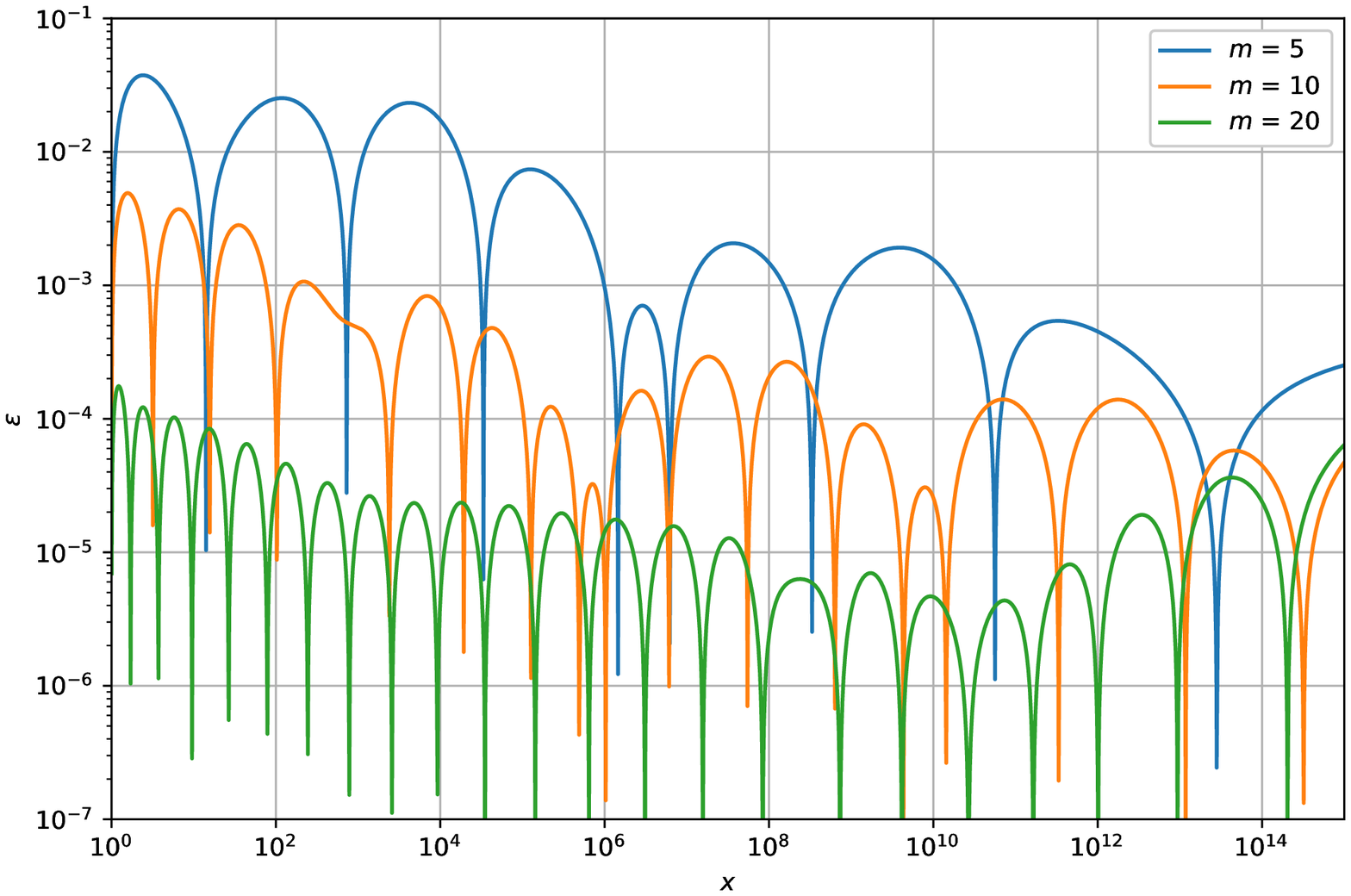} 
\includegraphics[width=0.49\linewidth]{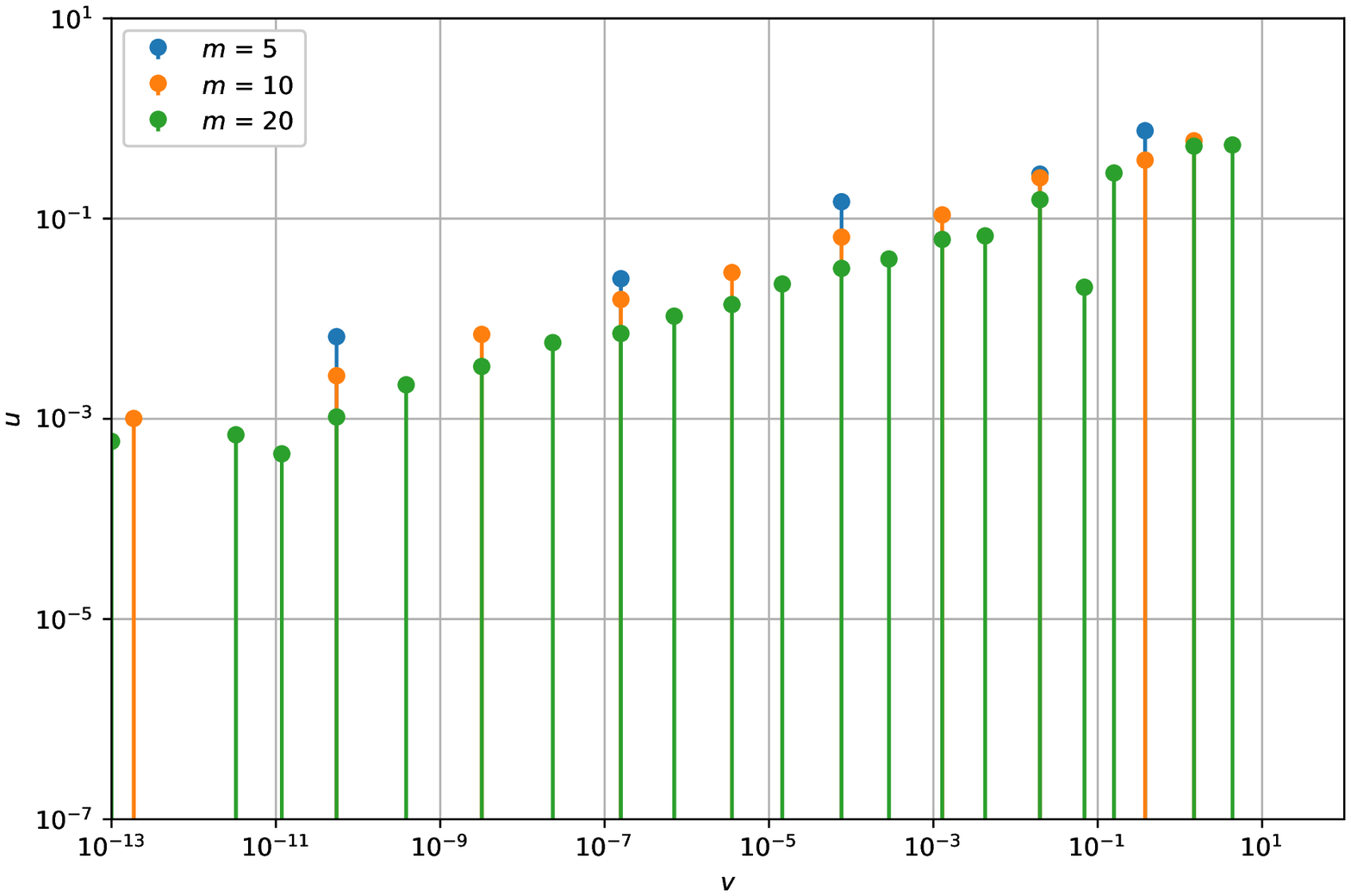} 
\caption{Approximation accuracy $\varepsilon$ (left) and approximation parameters $u_i, v_i, \ i = 1,2, \ldots, m,$ (right) for the function $x^{-\alpha}$ with $\alpha = 0.25$ at various $m$.} 
\label{fig-4}
\end{figure}

\begin{figure}[htbp]
\centering
\includegraphics[width=0.49\linewidth]{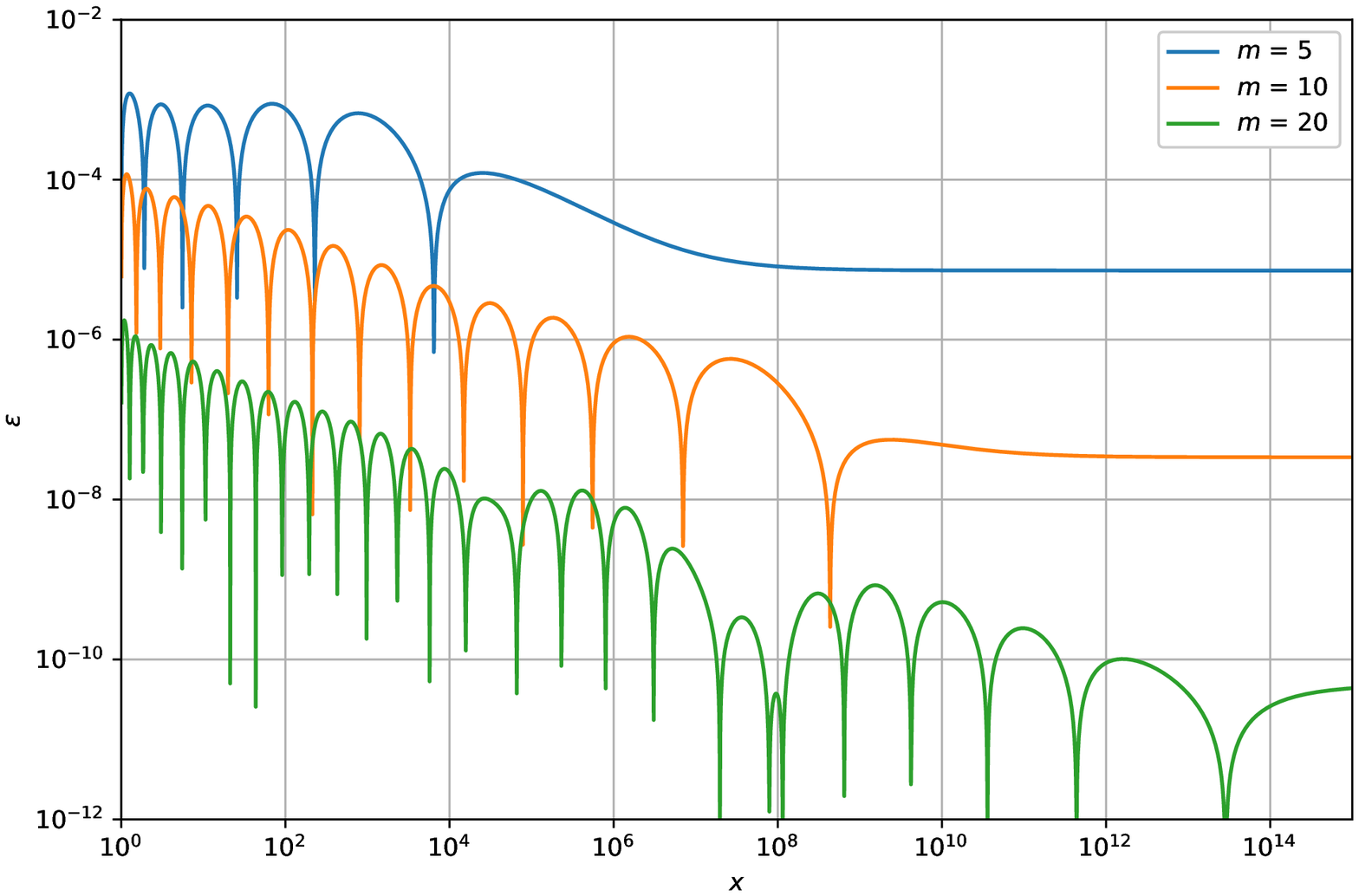} 
\includegraphics[width=0.49\linewidth]{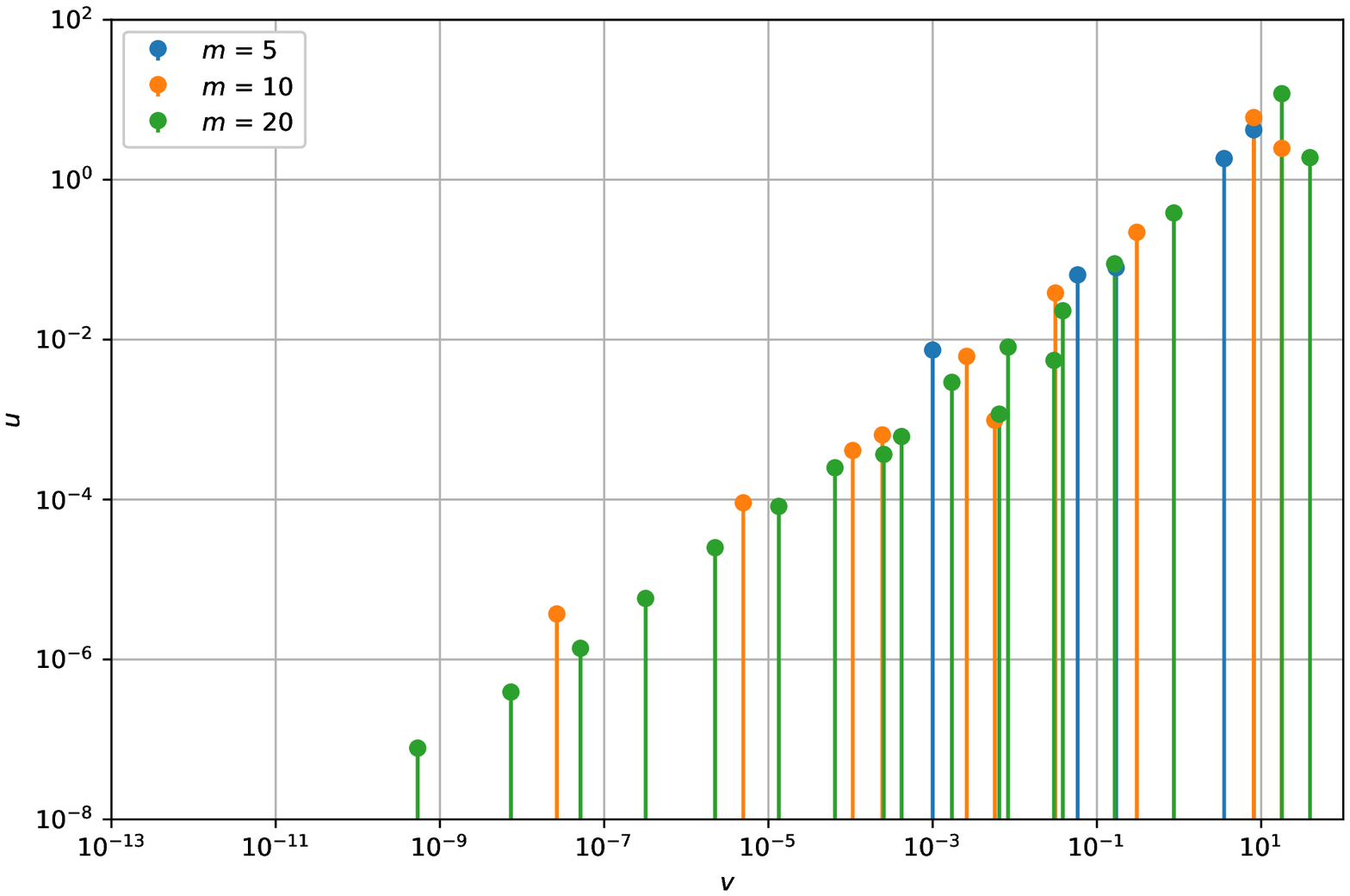} 
\caption{Approximation accuracy $\varepsilon$ (left) and approximation parameters $u_i, v_i, \ i = 1,2, \ldots, m,$ (right) for the function $x^{-\alpha}$ with $\alpha = 0.75$ at various $m$.} 
\label{fig-5}
\end{figure}

When approximating the function $x^{-\alpha}$, special attention is paid to the influence of the parameter $\alpha$.
The approximation accuracy for $m = 5, 10, 20$ at $\alpha = 0.25$ is shown in Figure \ref{fig-4}.
Similar results at $\alpha = 0.75$ are presented in Figure \ref{fig-5}.
As with other computational algorithms of rational approximation \cite{bonito2018numerical,harizanov2020rev}, increasing $\alpha$ increases the accuracy.
Note also that the calculation of the coefficient $u_i, v_i, \ i = 1,2, \ldots, m,$ at larger $\alpha$ is performed with a more significant number of iterations of the non-negative least squares method.
Table 1 shows the calculated values of the approximation parameter $u_i, v_i, \ i = 1,2, \ldots, 10,$ ($m = 10$) when approximating the function $x^{-\alpha}$ for $\alpha = 0.25, 0.5, 0.75$.

\begin{center}
\begin{table}[htp]
\label{tab-1}
\caption{Parameters approximation with $m=10$ for $x^{-\alpha}$ }
\centering
\begin{tabular}{l c ll c ll c ll}
\hline
 && \multicolumn{2}{l}{$\alpha  = $ 0.25} && \multicolumn{2}{l}{$\alpha  = $ 0.5}  && \multicolumn{2}{l}{$\alpha  = $ 0.75}\\
\cline{3-4} \cline{6-7} \cline{9-10} 
    $i$  && $u_i$  & $v_i$  &&  $u_i$   &  $v_i$    &&  $u_i$    &  $v_i$ \\
\hline
       1         &&  1.060084e-03     &  2.115485e-13     &&    1.263660e-04    &  5.816049e-09    &&  1.653295e-06    &  1.135126e-08  \\
       2         &&  2.778250e-03     &  6.526663e-11     &&    1.318851e-04    &  6.336196e-08    &&  1.664949e-05    &  6.273950e-07  \\
       3         &&  7.184790e-03     &  3.607348e-09     &&    2.177478e-03    &  2.389865e-06    &&  2.008706e-04    &  1.954833e-05  \\
       4         &&  1.608844e-02     &  1.812161e-07     &&    1.423375e-02    &  1.453310e-04    &&  9.792299e-04    &  2.343140e-04  \\
       5         &&  2.879614e-02     &  3.853128e-06     &&    3.605113e-02    &  3.090116e-03    &&  7.011612e-03    &  2.808580e-03  \\
       6         &&  6.751752e-02     &  8.192757e-05     &&    4.002657e-02    &  9.723689e-03    &&  3.444878e-02    &  3.059759e-02  \\
       7         &&  1.117978e-01     &  1.439033e-03     &&    7.561481e-02    &  4.933185e-02    &&  9.142663e-03    &  6.570394e-02  \\
       8         &&  2.518764e-01     &  2.088023e-02     &&    2.411886e-01    &  1.552328e-01    &&  2.280614e-01    &  3.333403e-01  \\
       9         &&  3.723954e-01     &  3.667548e-01     &&    3.672604e-01    &  1.397038e+00    &&  6.238136e+00    &  8.579865e+00  \\
      10         &&  6.229275e-01     &  1.537079e+00     &&    2.193928e+00    &  3.631519e+00    &&  2.478274e+00    &  1.842403e+01  \\
\hline
\end{tabular}
\end{table}
\end{center}

\clearpage

\subsection{Approximation of $\exp(- x^{\alpha} )$} 

The second example concerns the approximation of the stretching exponential function:
\[
f(x) = \exp(- x^{\alpha} ),
\quad 0 < \alpha < 1 .
\]
We will consider the case where the approximation is made with $a = 0, \ b = 10^3$, and an additional constraint (\ref{13}).
When approximating by the sum of exponents, we have
\[
\exp(- t^{\alpha}) \approx 1 + \sum_{i=1}^{m} u_i \varphi(x, v_i) ,
\]
where
\[
\varphi(x, v_i) = \exp( - v_i x) - 1,
\quad i = 1,2, \ldots, m .
\]

The choice of the weight function $\varrho(x)$ is done similarly to how it was when approximating the function $x^{-\alpha}$.
A uniform grid on the new variable $\theta \in [0, \beta]$ is introduced, with $x = \exp(\theta)-1$ ($\exp(\beta) - 1 = b$).
The calculations are performed using $\varrho(x) = (1+x)^{-1}$.

\begin{figure}[htbp]
\centering
\includegraphics[width=0.49\linewidth]{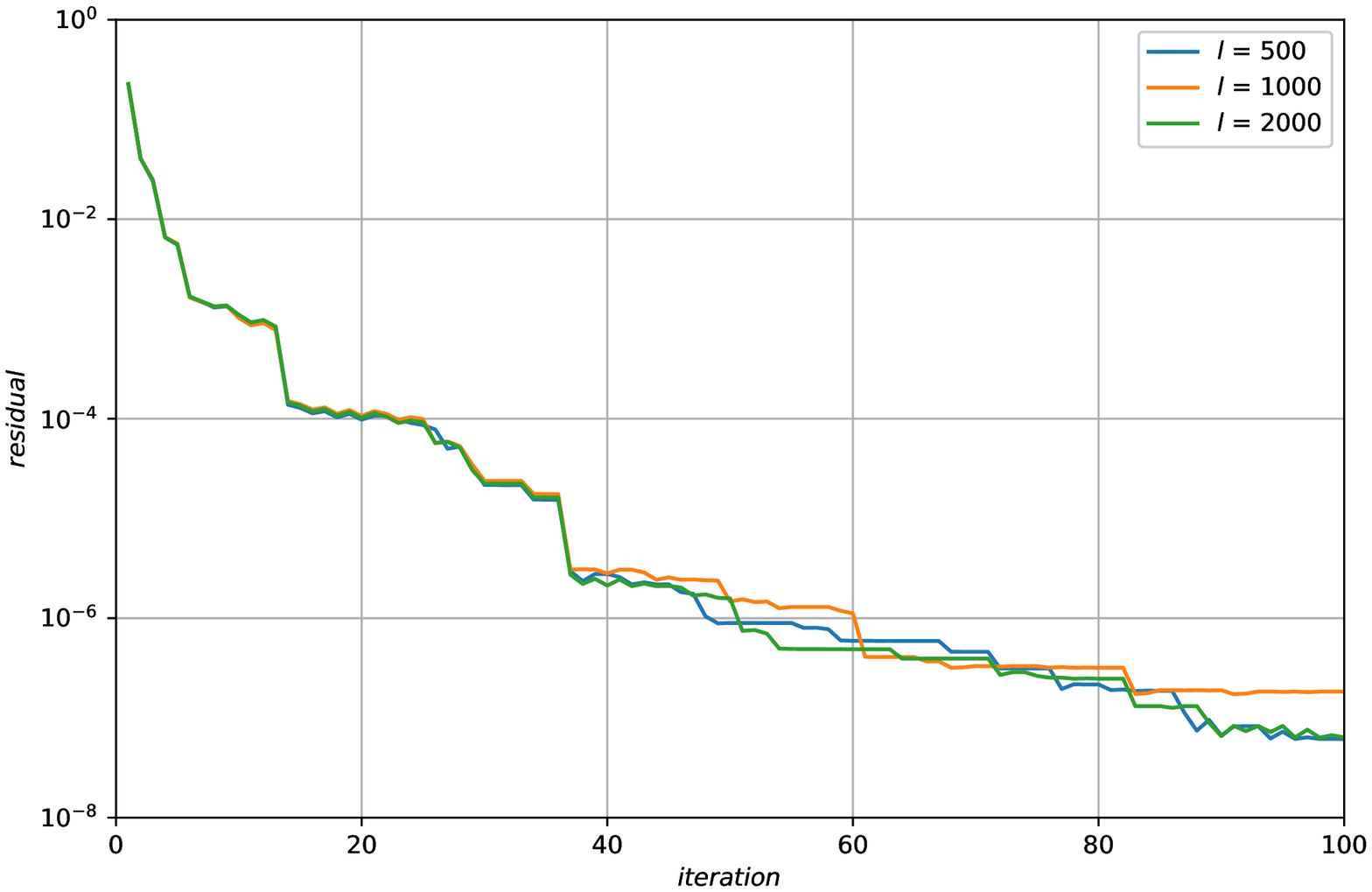} 
\includegraphics[width=0.49\linewidth]{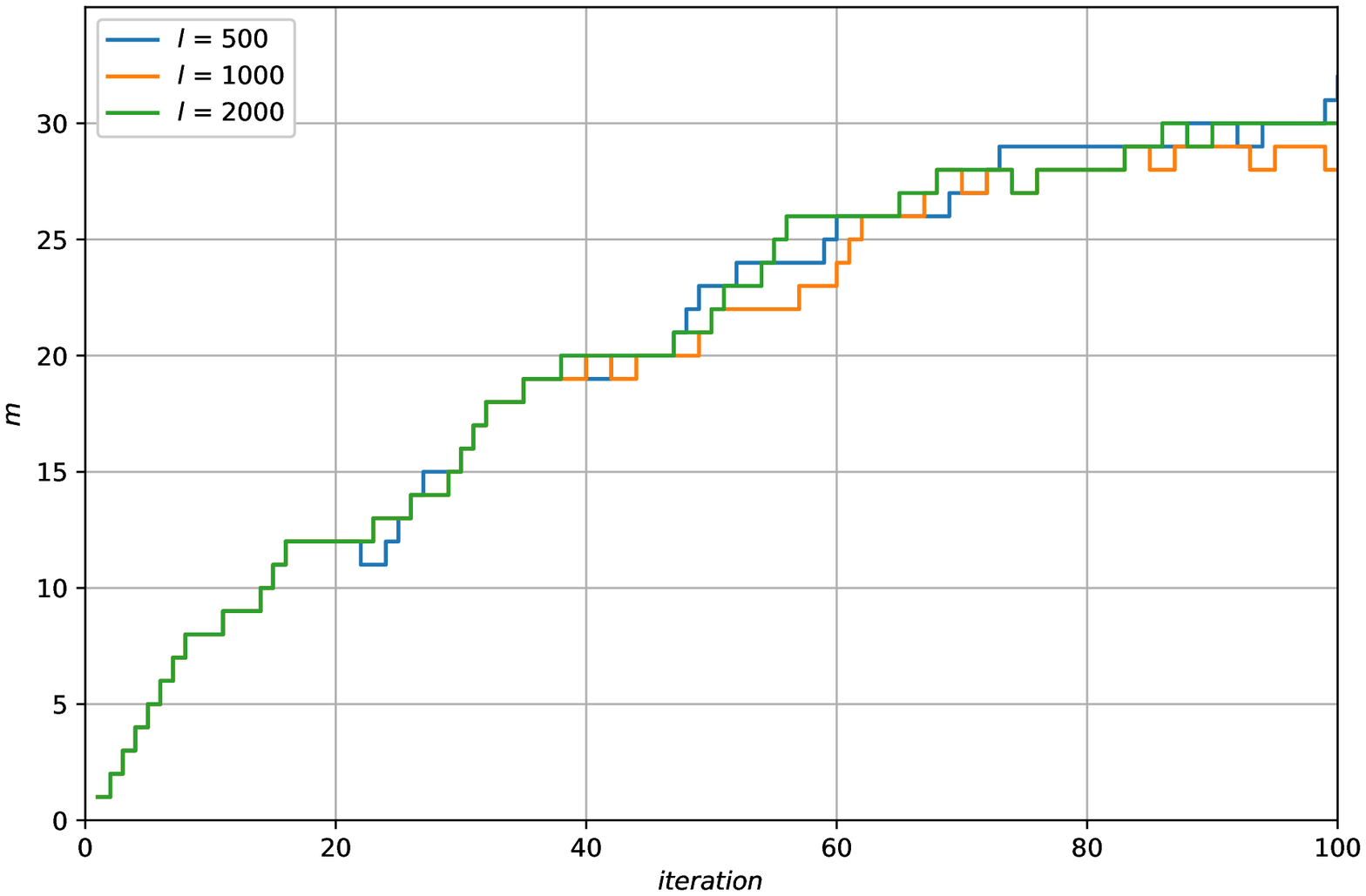} 
\caption{Residual (left) and number of nonzero elements $m$ (right) for $\widetilde{u}_k, \ k = 1,2, \ldots, l,$ in separate iterations of the NNLS method when approximating the function $\exp(- x^{\alpha} )$.}
\label{fig-6}
\end{figure}

\begin{figure}[htbp]
\centering
\includegraphics[width=0.49\linewidth]{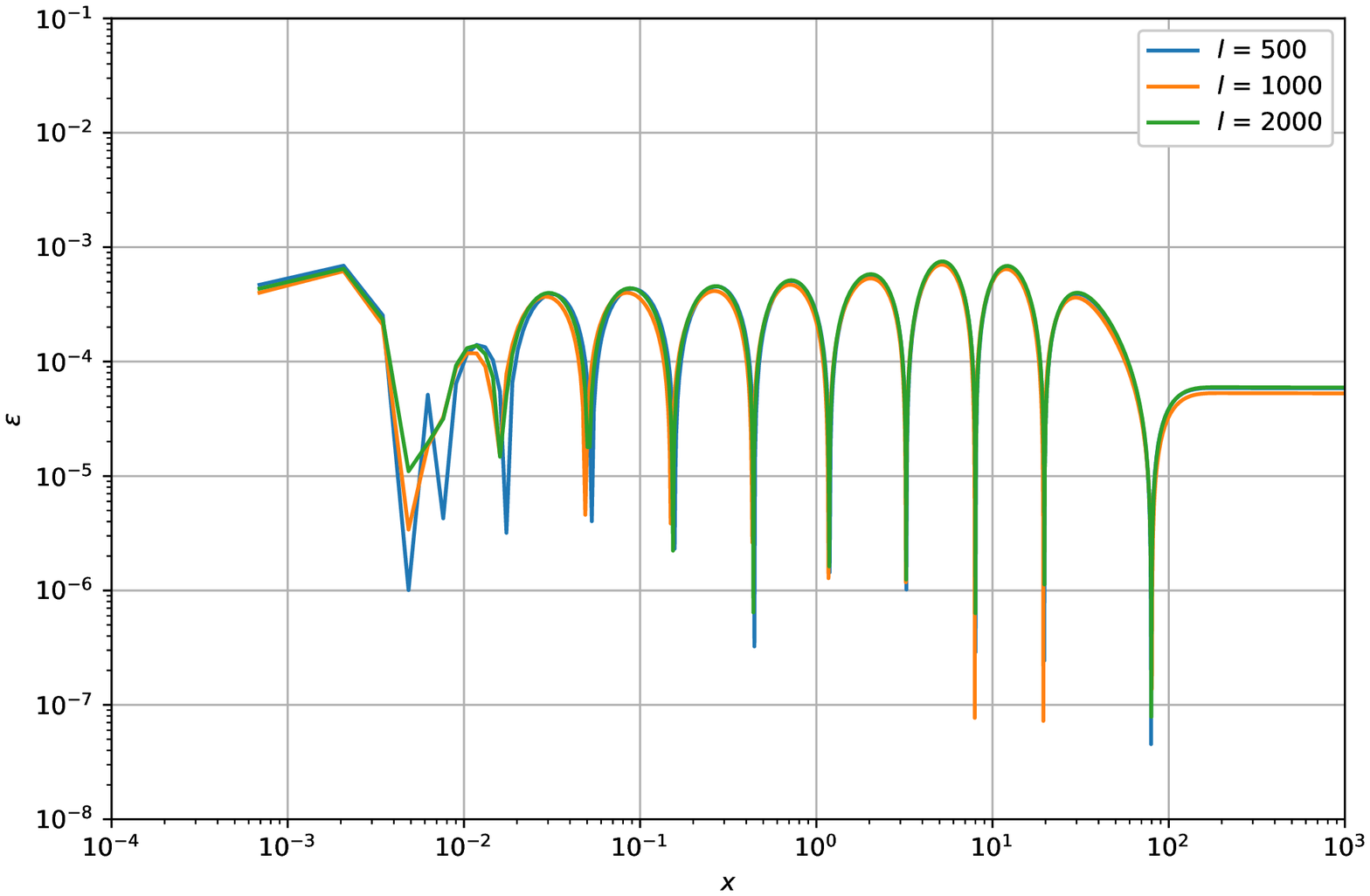} 
\includegraphics[width=0.49\linewidth]{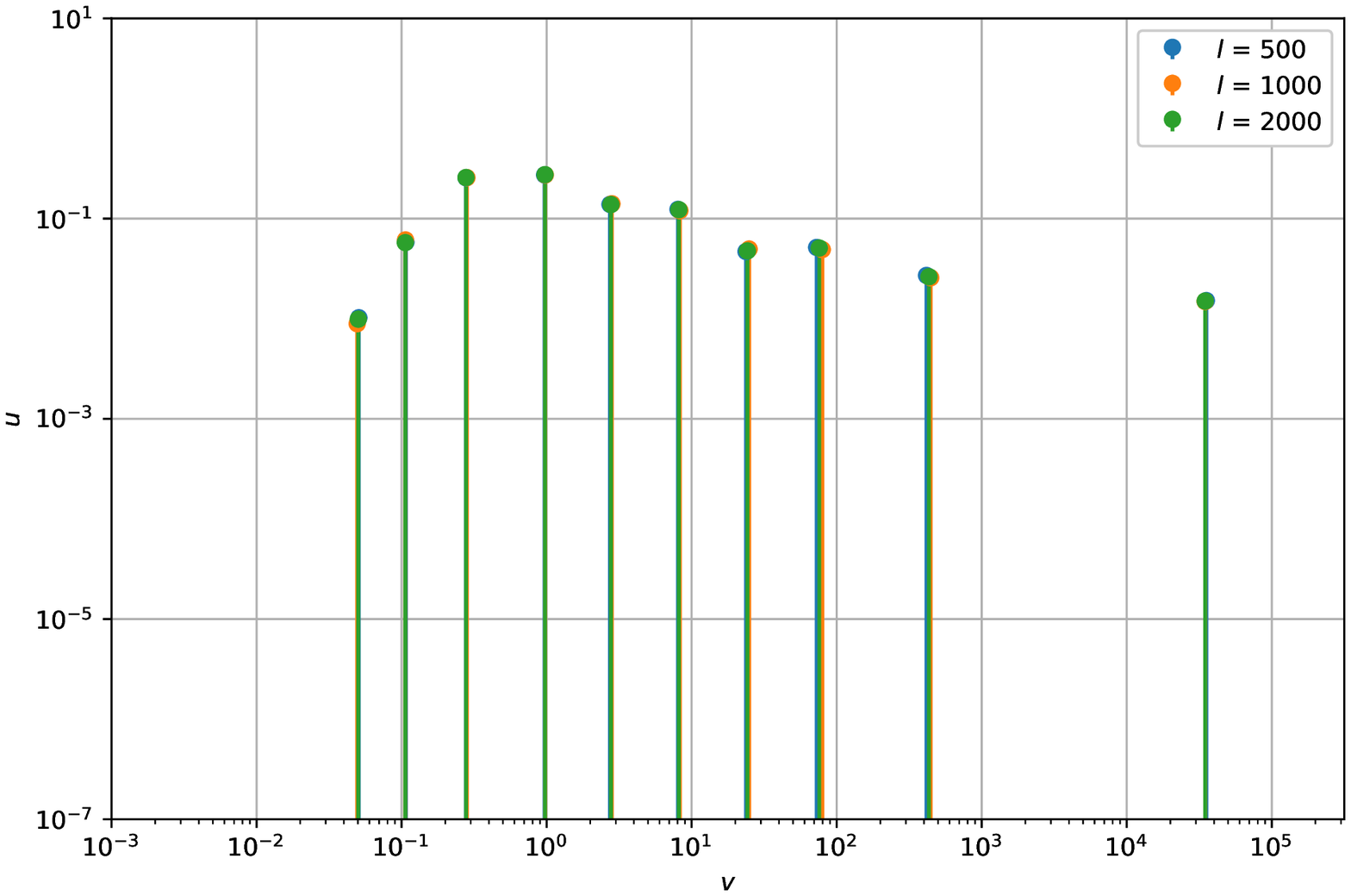} 
\caption{Approximation accuracy $\varepsilon$ (left) and approximation parameters $u_i, v_i, \ i = 1,2, \ldots, m,$ (right) at different partitioning of the interval $[c,d]$ when approximating $\exp(- x^{\alpha} )$.}
\label{fig-7}
\end{figure}

For the base variant the following parameters are chosen: $\alpha = 0.5$, $c = 10^{-4}$, $d = 10^{4}$, $n = 5000$, $l = 1000$, $m = 10$.
Figure~\ref{fig-6} shows the residual and the number of non-zero elements $m$ when using different numbers of iterations of the NNLS method.
The no significant influence of the interval partition detail of $[c, d]$ on the approximation accuracy and values of approximation parameters is illustrated by Fig.\ref{fig-7}.

\begin{figure}[htbp]
\centering
\includegraphics[width=0.49\linewidth]{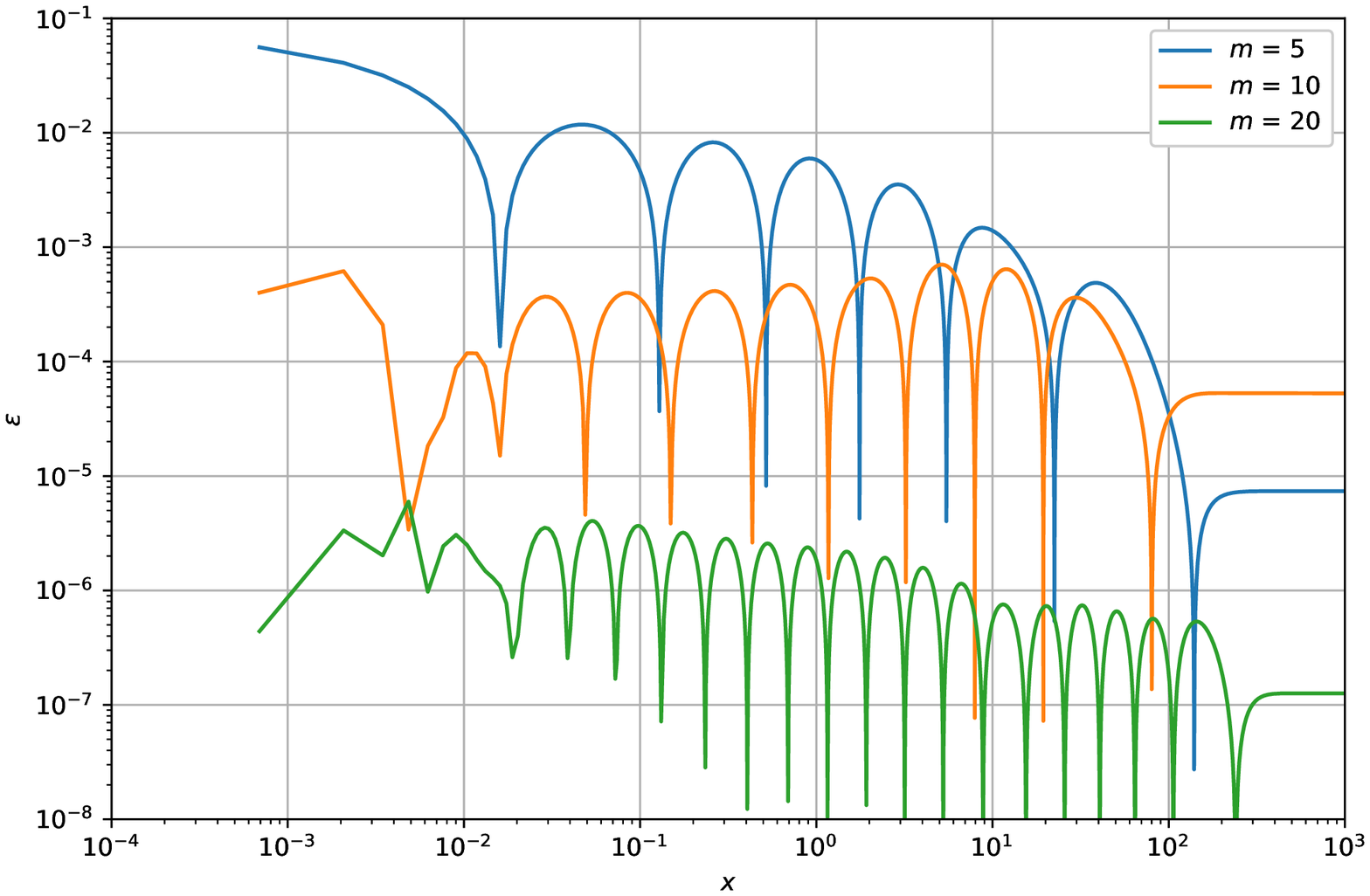} 
\includegraphics[width=0.49\linewidth]{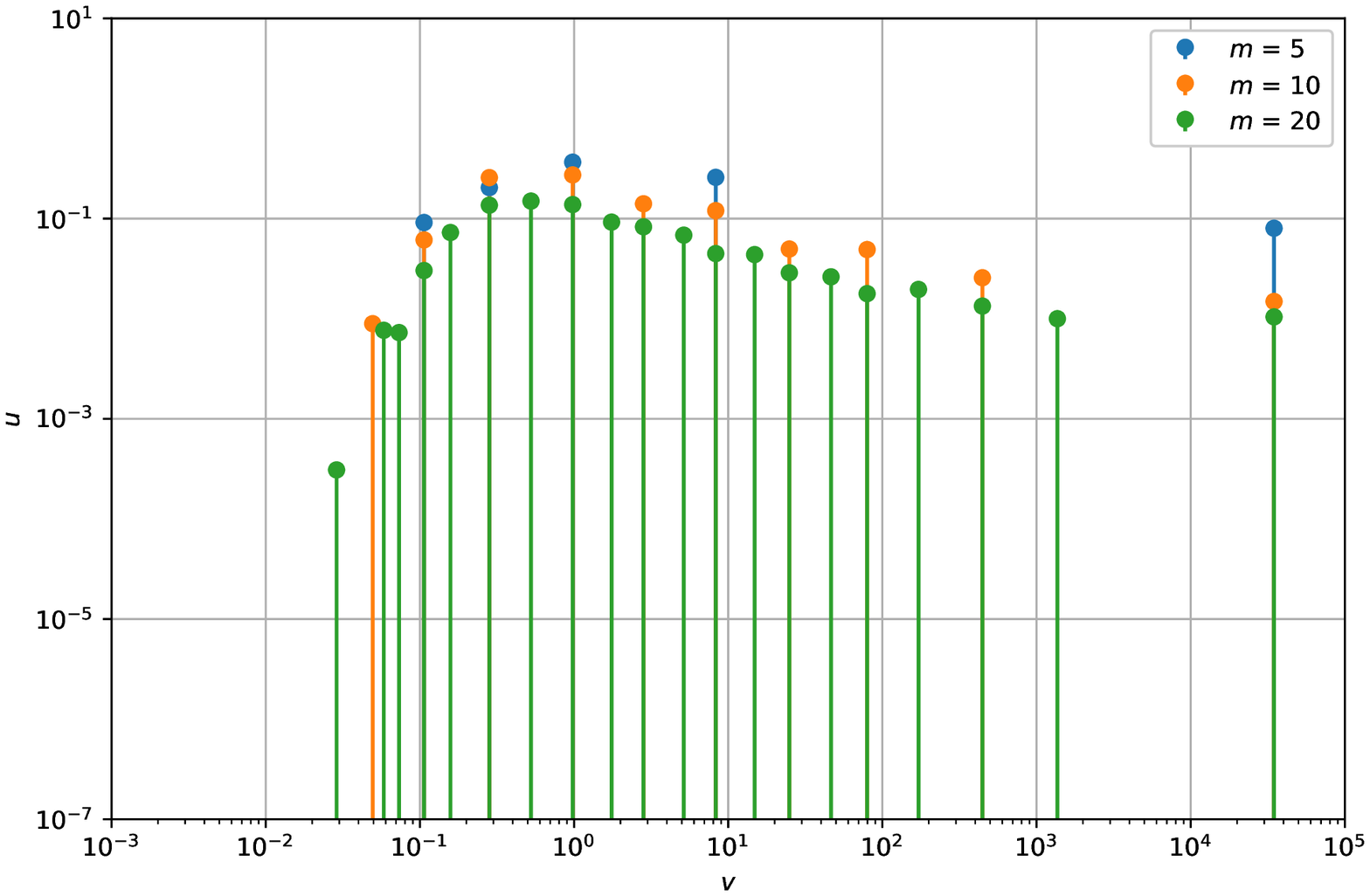} 
\caption{Approximation accuracy $\varepsilon$ (left) and approximation parameters $u_i, v_i, \ i = 1,2, \ldots, m,$ (right) for the function $\exp(- x^{\alpha} )$ with $\alpha = 0.5$ at various $m$.} 
\label{fig-8}
\end{figure}

\begin{figure}[htbp]
\centering
\includegraphics[width=0.49\linewidth]{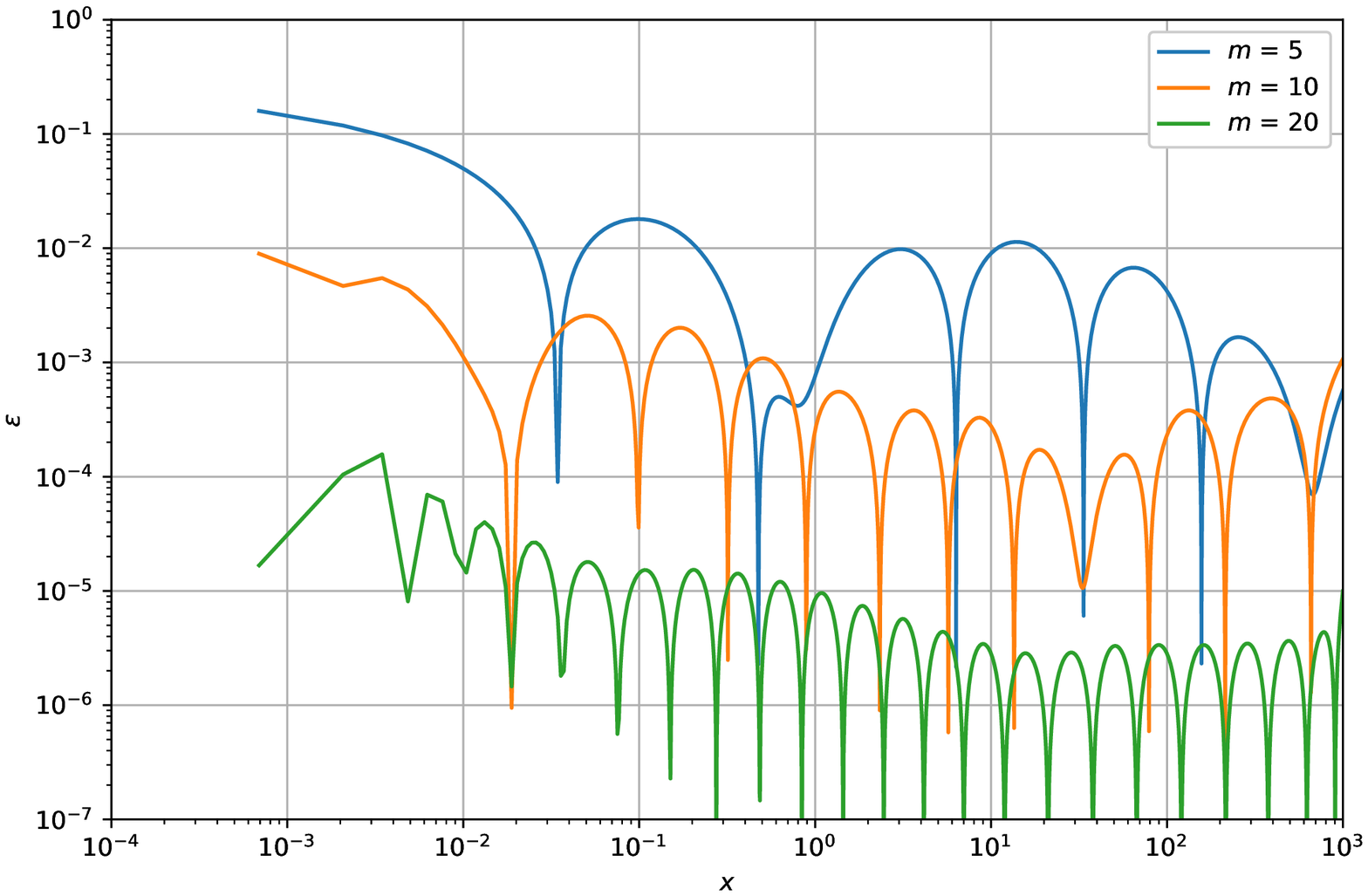} 
\includegraphics[width=0.49\linewidth]{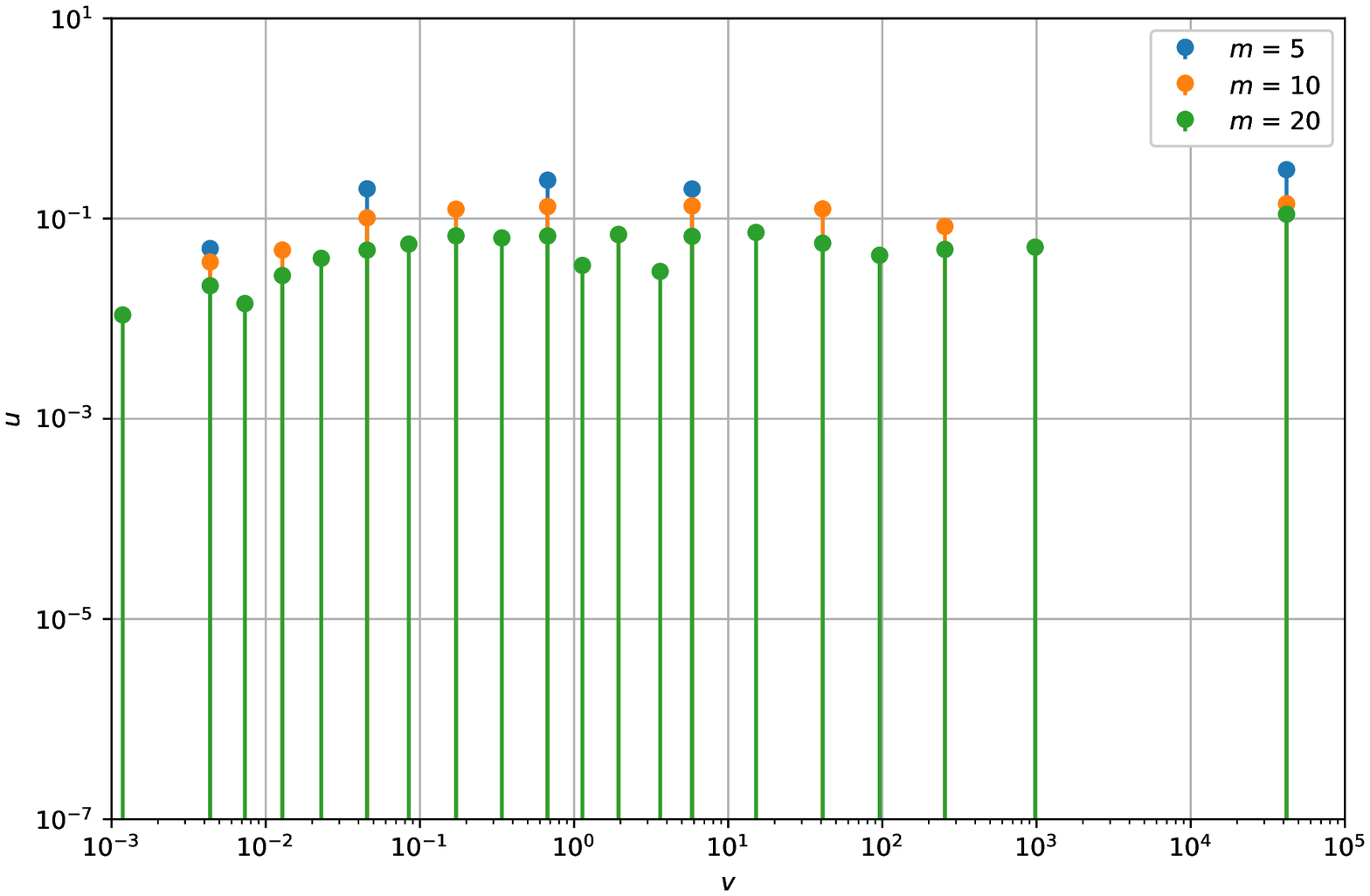} 
\caption{Approximation accuracy $\varepsilon$ (left) and approximation parameters $u_i, v_i, \ i = 1,2, \ldots, m,$ (right) for the function $\exp(- x^{\alpha} )$ with $\alpha = 0.25$ at various $m$.} 
\label{fig-9}
\end{figure}

\begin{figure}[htbp]
\centering
\includegraphics[width=0.49\linewidth]{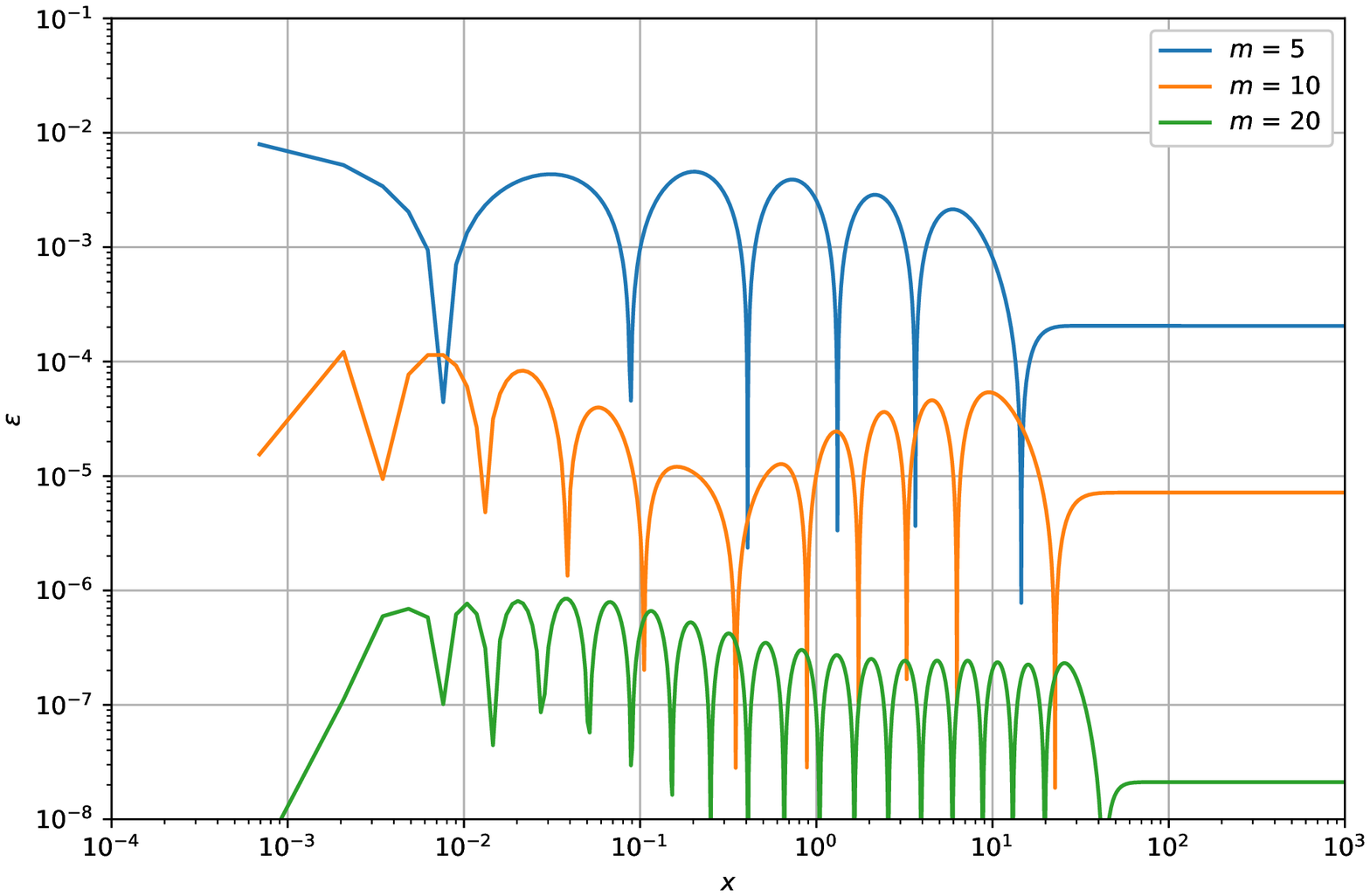} 
\includegraphics[width=0.49\linewidth]{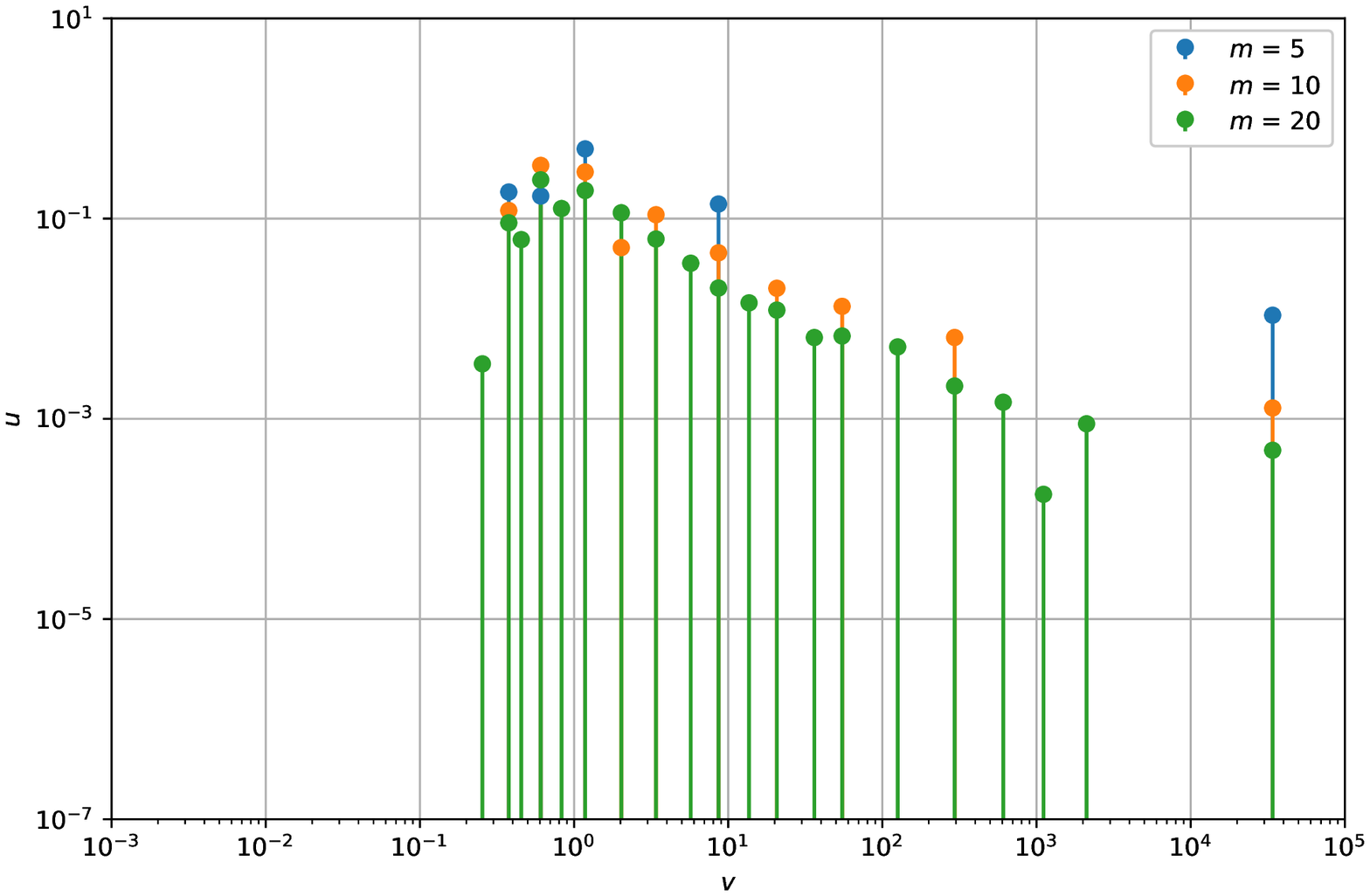} 
\caption{Approximation accuracy $\varepsilon$ (left) and approximation parameters $u_i, v_i, \ i = 1,2, \ldots, m,$ (right) for the function $\exp(- x^{\alpha} )$ with $\alpha = 0.75$ at various $m$.} 
\label{fig-10}
\end{figure}

The increase in approximation accuracy with increasing $m$ is illustrated by Fig.\ref{fig-8}; the figure also shows the approximation parameters
$u_i, v_i, \ i = 1,2, \ldots, m$. Similar data when setting $m = 5, 10, 20$ for the approximation of the function $\exp(- x^{\alpha} )$ at
$\alpha = 0.25$ and $\alpha = 0.75$ are shown in Figures \ref{fig-9},\ref{fig-10}, respectively.
When the parameter $\alpha$ is reduced, the approximation accuracy decreases.
The obtained approximation parameters $u_i, v_i, \ i = 1,2, \ldots, 10,$ ($m = 10$) for approximation of function $\exp(- x^{\alpha} )$ for $\alpha = 0.25, 0.5, 0.75$ are presented in Table 2.

\begin{center}
\begin{table}[htp]
\label{tab-2}
\caption{Parameters approximation with $m=10$ for $\exp(- x^{\alpha} )$}
\centering
\begin{tabular}{l c ll c ll c ll}
\hline
 && \multicolumn{2}{l}{$\alpha  = $ 0.25} && \multicolumn{2}{l}{$\alpha  = $ 0.5}  && \multicolumn{2}{l}{$\alpha  = $ 0.75}\\
\cline{3-4} \cline{6-7} \cline{9-10} 
    $i$  && $u_i$  & $v_i$  &&  $u_i$   &  $v_i$    &&  $u_i$    &  $v_i$ \\
\hline
       1         &&  3.684368e-02     &  4.361538e-03     &&    8.918599e-03    &  4.939622e-02    &&  1.204240e-01    &  3.772042e-01  \\
       2         &&  4.849511e-02     &  1.282650e-02     &&    6.127055e-02    &  1.064209e-01    &&  3.399225e-01    &  6.078323e-01  \\
       3         &&  1.017103e-01     &  4.546295e-02     &&    2.567263e-01    &  2.821308e-01    &&  2.922859e-01    &  1.180517e+00  \\
       4         &&  1.241971e-01     &  1.714882e-01     &&    2.731277e-01    &  9.794697e-01    &&  5.133306e-02    &  2.024447e+00  \\
       5         &&  1.319054e-01     &  6.742622e-01     &&    1.406392e-01    &  2.821308e+00    &&  1.093081e-01    &  3.400412e+00  \\
       6         &&  6.933268e-02     &  1.942175e+00     &&    1.200802e-01    &  8.296959e+00    &&  4.550720e-02    &  8.648423e+00  \\
       7         &&  1.337784e-01     &  5.831305e+00     &&    4.969377e-02    &  2.491130e+01    &&  2.014393e-02    &  2.066880e+01  \\
       8         &&  1.251127e-01     &  4.098384e+01     &&    4.906249e-02    &  7.959777e+01    &&  1.329286e-02    &  5.479472e+01  \\
       9         &&  8.343032e-02     &  2.543346e+02     &&    2.565896e-02    &  4.452959e+02    &&  6.492292e-03    &  2.940820e+02  \\
      10         &&  1.409798e-01     &  4.184289e+04     &&    1.487512e-02    &  3.471687e+04    &&  1.283043e-03    &  3.400412e+04  \\
\hline
\end{tabular}
\end{table}
\end{center}

\clearpage

\section{Conclusions}

\begin{enumerate}
\item The problem of nonlinear approximation of functions based on two sets of parameters is considered.
A particular term of the approximating function is the product of the first non-negative unknown parameter by a given nonlinear function that depends on the second unknown parameter.
\item  We present a heuristic computational algorithm for the nonlinear approximation of functions based on the minimization of the residual functional by values of the approximated function at separate points.
The unknown nonlinear approximation parameters are set on many points of the interval of permissible values.
Linear non-negative parameters are determined using the non-negative least squares method.
We find the solution of the nonlinear function approximation problem based on the estimation of the residual and the number of the first non-negative unknown parameters at each iteration of the non-negative least squares method.
\item The main points of practical use of the computational algorithm are illustrated by two typical problems of nonlinear approximation of functions. In the first example, we construct a rational approximation of the function $x^{-\alpha}, \ 0 < \alpha < 1$ at $x \geq 1$.
The second example is the approximation of the function $\exp(- x^{\alpha} ), \ \quad 0 < \alpha < 1$ at $ x \geq 0$ by the sum of exponents.
\end{enumerate}

\end{document}